\documentclass[12pt]{article}
\newcommand{\fracd}[2]{  \frac {\displaystyle {#1}}{\displaystyle {#2}}   }
\newcommand{\bb}{\begin{equation}}
\newcommand{\ee}{\end{equation}}
\newcommand{\esp}{\vspace*{0.5 cm}}

\begin{document}

\begin{center}
{\bf {\Large Recurrence Formulas for Fibonacci Sums}}
\end{center}

\begin{center}
Adilson J. V. Brand\~ao\footnote{Centro de Matem\'atica, Computa\c c\~ao e Cogni\c c\~ao,  Universidade Federal do ABC, Brazil. adilson.brandao@ufabc.edu.br}, 
Jo\~ao L. Martins \footnote{Departamento de Matem\'atica, Universidade Federal de Ouro Preto, 
Brazil. jmartins@iceb.ufop.br}
\end{center}

{\bf Abstract.} In this article we present a new recurrence formula for a finite sum involving the Fibonacci sequence. Furthermore, we state an algorithm to compute the sum of a power series related to Fibonacci series, without the use of term-by-term differentiation theorem 


\section{Introducion}
\renewcommand{\theequation}{1.\arabic{equation}}
\setcounter{equation}{0}

The Fibonacci sequence is one of the most famous numerical sequences in mathematics. It is defined in a recursive way: the first two terms are given and the following ones are defined as the sum of the two preceding ones. Mathematically speaking:
\[
F_0 = 0,~F_1 =1,~~F_r = F_{r-1} + F_{r-2},~r\geq 2  .
\]

The first terms are: $1, 1, 2, 3, 5, 8, 13, 21, \ldots $

This sequence comes from the single pair of rabbits' progeny problem, which was early proposed by Leonardo de Pisa (Fibonacci) at the Liber Abacci of 1202. An intriguing point is that this sequence appears in many problems from Mathematics as well as in Botanic, Crystallography, Computer Science, etc \cite{Dunlap}.

Consider the following finite sum involving the Fibonacci sequence, where $x$ is a real number, $m$ and $n$ are non-negative integers:
\bb \label{1.1}
\sum_{r=1}^{n} r^m F_r x^r.
\ee
Many authors have been seeking to establish a sum formula for (\ref{1.1}) 
(see \cite{Harris}, \cite{Brousseau}, \cite{Ledin},\cite{Gauthier}). In this article we state a sum formula for (\ref{1.1}) that we believe may be considered as a new result. 

Consider now the power series associated to (\ref{1.1}):
\bb \label{1.2}
\sum_{r=1}^{+\infty} r^m F_r x^r.
\ee
It is not difficult to demonstrate that the (\ref{1.2}) converges for all $m$ and all $x \in (-1/\phi , 1/\phi)$, in which $\phi=(1 + \sqrt5)/2$ is the golden ratio, a well-known constant associated with Fibonacci's sequence \cite{Dunlap}.

The question hereby interposed is the following: within its convergence interval, is there a formula for the sum of the (\ref{1.2}) series? An answer to this question is obtained by invoking the term-by-term differentiation theorem for power series.  Actually, such an equation is obtained by using $D = xd/dx$ operator $m$ times into the known identity
\[
\sum_{r=1}^{+\infty} F_r x^r = \fracd{x}{1 -x - x^2}.
\]
If we define
\[
S(x,j)=\sum_{r=1}^{+\infty} r^j F_r x^r ,
\]
a recurrence formula can be obtained by the following way: 
\bb \label{1.3}
\left \{
\begin{array}{ll}
\displaystyle{S(x,0)=\fracd{x}{1 -x - x^2}}, & \\
\displaystyle{S(x,j)=D[S(x,j-1)]}, & j=1, \ldots,m.
\end{array}
\right.
\ee

{\bf Example 1.1} Using the (\ref{1.3}) algorithm, we can calculate the numeric series' sum 
\bb \label{1.4}
\displaystyle{S=\sum_{r=1}^{+\infty} \fracd{r F_r}{3^r}.}
\ee
In fact, if  $S(x,0) = x/(1-x-x^2)$,  then $S(x,1) = x S' (x,0) = (x + x^3)/{(1 - x - x^2)}^2$. 
Hence, taking $x =1/3$ in $S(x,1)$,  we get the sum $S = 6/5$ for the (\ref{1.4}) series. 

{\bf Example 1.2} Try now to compute the numeric series' sum below by using the same algorithm:
\[
\displaystyle{\sum_{r=1}^{+\infty} \fracd{r^{50} F_r}{3^r}.}
\]
The (\ref{1.3}) algorithm's problem is, for each single step, higher computation cost in order to differentiate a function.  The example 1.2 points out this difficulty. In this paper we obtain another recurrence formula to calculate the sum of (\ref{1.2}).

The article is henceforth organized as follows: In the second section we present our main result, a recurrence formula for the (\ref{1.1}) finite sum and we show that it recovers some results on finite summation formulas involving the Fibonacci sequence. The third section was intended to rigorously proof our formula. In the fourth section we state an algorithm to compute the sum of (\ref{1.2}) without the use of derivatives. Finally, in the fifth section, we give some comments about the results and future possibilities.

\section{Finite Sums}
\renewcommand{\theequation}{2.\arabic{equation}}
\setcounter{equation}{0}
Our main result in this section is the theorem below:

{\bf Theorem 2.1.} Let $x \in \Re,~1-x- x^2 \not = 0$ be given. Then the following finite recorrence formula  holds

\bb \label{2.1}
\begin{array}{ll}
\displaystyle{\sum_{r=1}^{n} r^m F_r x^r = \fracd{1}{1 - x - x^2} \sum_{i=1}^m \pmatrix{m \cr i} {(-1)}^{i + 1} \sum_{r=1}^n r^{m - i} F_r x^r} + & \\
 & \\
+ \displaystyle{\fracd{x^2}{1 - x - x^2} \sum_{i=1}^m \pmatrix{m \cr i} \sum_{r=1}^{n-1} r^{m - i} F_r x^r - \fracd{n^m (F_{n+1} x^{n+1} + F_n x^{n+2})}{1 - x - x^2}} .& 
\end{array} 
\ee

\esp

As consequence of theorem 2.1 we obtain many closed formulas for finite sums involving Fibonacci summation. In fact,  taking $x=1$ in (\ref{2.1}) we obtain the following finite summation:
\bb \label{2.2}
\displaystyle{ \sum_{r=1}^{n} r^m F_r = \sum_{i=1}^m \pmatrix{m \cr i} {(-1)}^i \sum_{r=1}^n r^{m - i} F_r } -  \displaystyle{ \sum_{i=1}^m \pmatrix{m \cr i} \sum_{r=1}^{n-1} r^{m - i} F_r + n^m F_{n+2} }.  
\ee
We believe that (\ref{2.2}) is a new formula for (\ref{1.1}). From (\ref{2.2}) we can derive closed for some special cases of $m$. For instance, taking $m=1$ in 
(\ref{2.2}) we obtain
\bb \label{2.3}
\sum_{r=1}^{n} r F_r = \sum_{i=1}^1 \pmatrix{1 \cr i} {(-1)}^i \sum_{r=1}^n r^{1 - i} F_r  -   \sum_{i=1}^1 \pmatrix{1 \cr i} \sum_{r=1}^{n-1} r^{1 - i} F_r + n F_{n+2}.
\ee
It is well known (see \cite{Ledin}) that
\bb \label{2.4}
\displaystyle{ \sum_{r=1}^{n - 1} F_r} = F_{n+1} - 1. 
\ee
Thus, from (\ref{2.3}) and (\ref{2.4}) we conclude that
\[
\begin{array}{ll}
\displaystyle{\sum_{r=1}^{n} r F_r}  & =  - \displaystyle{\sum_{r=1}^{n}  F_r - \sum_{r=1}^{n - 1} F_r + n F_{n+2} } \\
 & \\
 & = -F_{n+2} + 1 - F_{n+1} + 1 + n F_{n+2}.
\end{array}
\]
Therefore
\bb \label{2.5}
\displaystyle{\sum_{r=1}^{n} r F_r = n F_{n+2} - F_{n+3} + 2}, 
\ee
which is the formula (1) that appears in \cite{Harris}. Now, taking $m=2$ in (\ref{2.2}) we can see  that
\[
\sum_{r=1}^{n}r^{2}F_{r}=\sum_{i=1}^{2}\pmatrix{2 \cr i}(-1)^{i}\sum_{r=1}^{n} \, 
\,r^{2-i}F_{r}
-\sum_{i=1}^{2}\pmatrix{2 \cr i}\sum_{r=1}^{n-1}\,r^{2-i}F_{r}
+n^{2}F_{n+2} ,
\]
that is,
\bb \label{2.6}
\sum_{r=1}^{n}r^{2}F_{r}=-2\sum_{r=1}^{n}rF_{r}-2\sum_{r=1}^{n-1}rF_{r} +n^{2}F_{n+2}+\sum_{r=1}^{n}F_{r}-\sum_{r=1}^{n-1}F_{r}.
\ee
Thus, using (\ref{2.4}) and (\ref{2.5}) in (\ref{2.6}), after some algebric manipulation, we obtain
\[
\sum_{r=1}^{n}r^2F_{r}=(n^2+2)F_{n+2}-(2n-3)F_{n+3}-8 ,
\]
which is the formula (17) in \cite{Harris}. In an analogous way we can recover other known identities taking different values for $m$ in (\ref{2.2}). Actually, the recorrence formula 
(\ref{2.1}) can produce a lot of identities, simply choosing special values to $x$ and $m$. For instance, the formula (\ref{2.1}) for $x=-1$ is  
\bb \label{2.7}
\begin{array}{lll}
\displaystyle{\sum_{r=1}^{n}(-1)^{r}r^{m}F_{r}}  & = &\displaystyle{\sum_{i=1}^{m}\pmatrix{m \cr i}(-1)^{i+1}\sum_{r=1}^{n} \, (-1)^{r}\,r^{m-i}F_{r} }+\nonumber \\
&+&\displaystyle{ \sum_{i=1}^{m}\pmatrix{m \cr i}\sum_{r=1}^{n-1}\,(-1)^{r}\,r^{m-i}F_{r}
-n^{m}\left[(-1)^{n+1}F_{m+1}+(-1)^{n+2}F_{n}\right ]}.
\end{array}
\ee
Taking $m=1$ in (\ref{2.7}) we obtain
\bb \label{2.8}
\sum_{r=1}^{n}(-1)^{r}\,r\,F_{r}=\sum_{r=1}^{n}(-1)^{r}F_{r}+\sum_{r=1}^{n-1}(-1)^{r}F_{r}
-n\left[(-1)^{n+1}F_{n+1}+(-1)^{n+2}F_{n}\right] .
\ee
Since that (see \cite{Harris})
\bb \label{2.9}
\sum_{r=1}^{n-1}(-1)^{r}F_{r} = (-1)^{n-1} F_{n-2} -1,
\ee
and using (\ref{2.9}) in (\ref{2.8}), we conclude that
\bb \label{2.10}
\sum_{r=1}^{n}(-1)^{r}rF_{r}=(-1)^{n}F_{n-1} -1 + (-1)^{n-1}F_{n-2} -1 - n\left[(-1)^{n+1}F_{n+1}+(-1)^{n+2}F_{n}\right].
\ee
After some simplifications, (\ref{2.10}) becomes
\bb \label{2.11}
\sum_{r=1}^{n}(-1)^{r}rF_{r}=(-1)^{n} (n+1) F_{n-1} + (-1)^{n-1} F_{n-2} -2.
\ee
Note that (\ref{2.11}) is the formula (2) in \cite{Harris}.

\section{Proof of Our Main Result}

\renewcommand{\theequation}{3.\arabic{equation}}
\setcounter{equation}{0}

Before proving Theorem 2.1 we need to state some auxiliar results: 

{\bf Lemma 3.1.}  Let non-negative integers $n\geq k$ be given. Suppose that $1 - x - x^2 \not =0$. Then
\bb \label{3.1}
F_k x^k + F_{k+1} x^{k+1} + \ldots + F_n x^n = 
\fracd{F_k x^k + F_{k-1} x^{k+1} - F_{n+1} x^{n+1} - F_n x^{n+2}}{1 - x - x^2} .
\ee
{\sl Proof.} Consider the sum
\bb \label{3.2}
S = F_k x^k + F_{k+1} x^{k+1} + F_{k+2}x^{k + 2} +  F_{k+3}x^{k+3} + \ldots  + F_{n-1}x^{n-1} + F_n x^n.
\ee
Multiplying (\ref{3.2}) by $-x$ and $-x^2$ we obtain
\bb \label{3.3}
-x S = -F_k x^{k +1} - F_{k+1} x^{k + 2} - F_{k+2}x^{k + 3}- \ldots - F_{n-2}x^{n - 1} - 
F_{n-1}x^{n } - F_n x^{n + 1}.
\ee
\bb \label{3.4}
-x^2 S = -F_k x^{k +2} - F_{k+1} x^{k+3} - \ldots - F_{n-3} x^{n-1}
 - F_{n-2}x^{n} - F_{n-1}x^{n + 1}- F_n x^{n +2}.  
\ee
Adding (\ref{3.2}), (\ref{3.3}) and (\ref{3.4}),  remembering the definition of Fibonacci sequence and cancelling terms we have
\bb
S -xS -x^2S =F_k x^k + F_{k+1} x^{k+1} -F_k x^{k +1} - F_n x^{n + 1}  - F_{n-1}x^{n + 1}- F_n x^{n +2}
\ee
Using again the definition of Fibonacci sequence we conclude that
\[
S -xS -x^2S =F_k x^k + F_{k-1} x^{k+1} - F_{n+1} x^{n+1} - F_n x^{n+2},
\]
that is, (\ref{3.1}) holds. 

{\bf Lemma 3.2.} Let $x \in \Re,~1-x- x^2 \not = 0$ be given. Then the following identity holds
\bb
\begin{array}{ll}
\displaystyle{\sum_{r=1}^n r^m F_r x^r} = & \displaystyle{\fracd{1}{1 - x - x^2} \sum_{r=1}^n ( r^m - {(r-1)}^m ) ( F_r x^r + F_{r-1} x^{r+1})} - \\
 & \\
 & - \displaystyle{\fracd{n^m (F_{n+1} x^{n+1} + F_n x^{n+2})}{1 - x - x^2}} .
\end{array}
\ee
{\sl Proof.} Consider the sum
\[
\sum_{r=1}^n r^m F_r x^r= 1 F_1 x + 2^m F_2 x^2 + 3^m F_3 x^3 + \ldots + n^m F_n x^n .
\]
It is easy to see that the sum above can be rearranged in the following way
\[
\sum_{r=1}^n r^m F_r x^r = (F_1 x +  F_2 x^2 + F_3 x^3 + \ldots +  F_n x^n) +
\]
\[
+ (2^m - 1)(F_2 x^2 + F_3 x^3 + \ldots +  F_n x^n) +
\]
\[
+ (3^m - 2^m)(F_3 x^3 + \ldots +  F_n x^n) + \ldots +
\]
\[
+ ((n-1)^m - (n-2)^m)(F_{n-1} x^{n-1}+  F_n x^n) +
\]
\[
+ (n^m - (n-1)^m)( F_n x^n).
\]
By using the Lemma 3.1 we can write the last sum as
\[
\sum_{r=1}^n r^m F_r x^r = \fracd{F_1 x^1 + F_0 x^2 - F_{n+1} x^{n+1} - F_n x^{n+2}}{1 - x - x^2} +
\]
\[
+ (2^m - 1)(\fracd{F_2 x^2 + F_1 x^3 - F_{n+1} x^{n+1} - F_n x^{n+2}}{1 - x - x^2}) + \ldots +
\]
\[
+ (n^m - (n-1)^m)(\fracd{F_n x^n + F_{n-1} x^{n+1} - F_{n+1} x^{n+1} - F_n x^{n+2}}{1 - x - x^2}).
\]
Therefore, the sum can be expressed by 
\[
\sum_{r=1}^n r^m F_r x^r = \fracd{1}{1 - x - x^2} \sum_{r=1}^n ( r^m - {(r-1)}^m ) ( F_r x^r + F_{r-1} x^{r+1}) -
\]
\[
- [\fracd{F_{n+1} x^{n+1} + F_n x^{n+2}}{1 - x - x^2}] [ 1 + (2^m -1) + 
(3^m - 2^m) + \ldots + ( {(n-1)}^m - {(n-2)}^m ) +  ( n^m - {(n-1)}^m )  ]
\]
After cancelling some terms we finally obtain
\[
\sum_{r=1}^n r^m F_r x^r = \fracd{1}{1 - x - x^2} \sum_{r=1}^n ( r^m - {(r-1)}^m ) ( F_r x^r + F_{r-1} x^{r+1}) -
\]
\[
- \fracd{n^m (F_{n+1} x^{n+1} + F_n x^{n+2})}{1 - x - x^2},
\]
which is the desired result.

{\sl Proof of Theorem 2.1.} By using a suitable change of variables we have 
\[
\fracd{1}{1 - x - x^2}  \sum_{r=1}^n ( r^m - {(r-1)}^m ) ( F_r x^r + F_{r-1} x^{r+1}) 
\]
\[
=\fracd{1}{1 - x - x^2} \sum_{r=1}^n ( r^m - {(r-1)}^m ) F_r x^r
+ \fracd{1}{1 - x - x^2} \sum_{\theta=1}^n ( \theta^m - {(\theta-1)}^m )F_{\theta-1} x^{\theta+1}) 
\]
\[
=\fracd{1}{1 - x - x^2}\sum_{r=1}^n ( r^m - {(r-1)}^m ) F_r x^r
+ \fracd{x^2}{1 - x - x^2} \sum_{r=1}^{n - 1} ( {(r+1)}^m - r^m ) F_r x^r 
\]
\[
= \fracd{1}{1 - x - x^2}\sum_{r=1}^n \sum_{i=1}^m \pmatrix{m \cr i} {(-1)}^{i + 1} r^{m - i} F_r x^r + \fracd{x^2}{1 - x - x^2}\sum_{r=1}^{n-1} \sum_{i=1}^m \pmatrix{m \cr i} r^{m - i} F_r x^r 
\]
\[
= \fracd{1}{1 - x - x^2} \sum_{i=1}^m \pmatrix{m \cr i} {(-1)}^{i + 1} \sum_{r=1}^n r^{m - i} F_r x^r + \fracd{x^2}{1 - x - x^2} \sum_{i=1}^m \pmatrix{m \cr i} \sum_{r=1}^{n-1} r^{m - i} F_r x^r 
\]
The theorem follows by the result above and the Lemma 3.2.

\section{Power Series} 
\renewcommand{\theequation}{4.\arabic{equation}}
\setcounter{equation}{0}

In this section we state a result which provides an algorithm to compute the sum of 
(\ref{1.2}) without the use of term-by-term differentiation theorem. This algorithm is a consequence of the theorem below:

{\bf Theorem 4.1.} Let $x \in (-1/ \phi,1/ \phi)$ be given. Then the following recorrence formula  holds

\bb \label{4.1}
\begin{array}{ll}
\displaystyle{\sum_{r=1}^{+\infty} r^m F_r x^r} = & \displaystyle{\fracd{1}{1 - x - x^2} 
\sum_{i=1}^m \pmatrix{m \cr i} {(-1)}^{i + 1} \sum_{r=1}^{+\infty} r^{m - i} F_r x^r} +  \\
 &  \\
 & + \displaystyle{\fracd{x^2}{1 - x - x^2}  \sum_{i=1}^m \pmatrix{m \cr i} \sum_{r=1}^{+\infty} r^{m - i} F_r x^r } .
\end{array} 
\ee
{\sl Proof.} Considering the Theorem 2.1, it is sufficient to take $n \rightarrow +\infty$ in (\ref{2.1}) and to remember that $\displaystyle{\lim_{n \rightarrow +\infty} n^m F_n x^n = 0}$, since the series (\ref{1.2}) converges for all integer $m$ and $x \in (-1/ \phi,1/ \phi)$. 

\esp

By theorem 4.1, we can obtain the following algorithm in order to provide the sum of 
(\ref{1.2}):
\bb \label{4.2}
\left \{
\begin{array}{l}
S(x,0)=  \fracd{x}{1 -x - x^2},  \\
S(x,j)=  \displaystyle{\fracd{1}{1 - x - x^2} \sum_{i=1}^j \pmatrix{j \cr i} {(-1)}^{i + 1} S(x,j - i)}+ \displaystyle{\fracd{x^2}{1 - x - x^2}  \sum_{i=1}^j \pmatrix{j \cr i} S(x,j - i)}, \\
j=1, \ldots,m.    
\end{array}
\right.
\ee

This algorithm can be implemented in an efficient way, instead of the expensive process using the standard derivative operator. It answers, for instance, the question proposed in the example 1.2: 
\[
\displaystyle{\sum_{n=1}^{+\infty}\fracd{n^{50} F_n}{3^n} = 6.526\times 10^{74}}.
\]

\section{Final Remarks}
\renewcommand{\theequation}{5.\arabic{equation}}
\setcounter{equation}{0}

In this article we state a new recurrence formula for a finite sum related to Fibonacci sequence. This formula recovers a lot of identities for Fibonacci sums. Besides this, it implies an algorithm to compute the sum of Fibonacci power series without the use of derivatives. The scheme used to obtain this results can be extended to others  series. The ideas presented here are part of a larger investigation which has been developed concerning the series
\bb \label{5.1}
\displaystyle{\sum_{r=1}^{+\infty}}{r^{m}}x^{r}a_{r},
\ee 
in which $\{a_r\}$ is an arbitrary sequence. In this article $\{a_r\}$ is the Fibonacci sequence.  Nevertheless, we can extend our results for other sequence types (see \cite{JB1}, \cite{JB2}). For example, if we take $a_r  = 1$, (\ref{5.1}) turns into the generalized geometric series 
\bb
\sum_{r=1}^{+\infty} r^{m} x^r,
\ee 
which converges for all $ x \in (-1,1)$.  Using the same ideas developed in the last section, we can find out a recurrence formula for such a series:
\[
\sum_{r=1}^{+\infty} r^{m} x^{r} = \fracd{1}{1 - x} \sum_{i=1}^m \pmatrix{m \cr i} {(-1)}^{i + 1} 
\sum_{r=1}^{+\infty} r^{m - i} x^r \cdot   
\]
There are other subjects still under investigation by which we search to extend the results hereby presented for other sequences such as, Lucas', Generalized Fibonacci's, Generalized Lucas',  Pell's, Tribonacci's sequences, etc. It should be observed that in 
\cite{Filipponi}, the author studied a series related to (\ref{1.2}), covering Lucas' and Fibonacci's generalized sequences. However, their results are only valid for a positive rational $x$. Besides, the employed technique is quite different from ours.  

Additional references concerning Fibonacci numbers and the golden ratio can be found in \cite{Dunlap}.

\end{document}